\documentclass[11pt]{article}

\usepackage{a4wide,authblk,epsfig,amsmath,amssymb,makeidx,graphics,color,pgfplots,stmaryrd,tikz}              
\usepackage{color}                   
\usetikzlibrary{calc}           
\usetikzlibrary{patterns}

\newtheorem{Th}{Theorem}[section]
\newtheorem{Le}[Th]{Lemma}
\newtheorem{Co}[Th]{Corollary} 
\newtheorem{Pro}[Th]{Proposition}

\newtheorem{Def}[Th]{Definition}
\newtheorem{rem}[Th]{Remark}

\newcommand{\be}{\begin{equation}}
\newcommand{\ee}{\end{equation}}
\newcommand{\hdrie}{\hspace{3mm}}
\newcommand{\und}{\hdrie\mbox{\rm and }\hdrie}

\author{Jan Brandts and Michal K\v{r}\'{\i}\v{z}ek\\
In memory of Professor Miroslav Fiedler (1926-2015)}  
\title{Factorization of CP-rank-$3$ Completely Positive Matrices}
\date{\today} 
       
\newcommand{\RR}{\mathbb{R}}
\newcommand{\CC}{\mathcal{C}} 
\newcommand{\TT}{\mathbb{T}}
\newcommand{\PP}{{\mathcal P}}
\newcommand{\sphere}{\mathbb{S}} 
\newcommand{\UU}{\mathcal{U}}
\newcommand{\FF}{\mathcal{F}}
\newcommand{\GG}{\mathcal{G}}
\newcommand{\Rdp}{\RR^3_{\geq 0}}
\newcommand{\st}{\sqrt{2}}
\newcommand{\sth}{\hspace{2mm} | \hspace{2mm}}
\newcommand{\half}{\frac{1}{2}}
\newcommand{\Rnspd}{\RR^{n\times n}_{\rm spd}}
\newcommand{\DD}{\mathcal{D}^{n\times n}}

\begin{document}
\maketitle 

\vspace*{-7mm} 
\begin{abstract}
A symmetric positive semi-definite matrix $A$ is called completely positive if there exists a matrix $B$ with nonnegative entries such that $A=BB^\top$. If $B$ is such a matrix with a minimal number $p$ of columns, then $p$ is called the cp-rank of $A$. In this paper we develop a finite and exact algorithm to factorize any matrix $A$ of cp-rank $3$. Failure of this algorithm implies that $A$ does not have cp-rank $3$. Our motivation stems from the question if there exist three nonnegative polynomials of degree at most four that vanish at the boundary of an interval and are orthonormal with respect to a certain inner product.
\end{abstract}
 
\section{Introduction and motivation} 
We study the problem of isometric embedding of a finite point set in $\RR^3$ into the nonnegative octant $\RR^3_{\geq 0}$. This problem appears in a number of applications, among which the factorization of {\em completely positive matrices} \cite{BeSh}. These are matrices $A$ that can be written as $A=BB^{\top}$ for some matrix $B$ with nonnegative entries. In Section \ref{sect1.2} we comment on the origin of our interest in this problem, which is connected to the existence of orthonormal bases of nonnegative functions, for instance, polynomials. First we review some known facts \cite{BeSh}.
 
\subsection{Isometric embedding of vectors into the nonnegative orthant}\label{sect-1.2}
We will write $\Rnspd$ for the set of real symmetric positive semidefinite $n\times n$ matrices, and 
\[ \DD=\{A\in\Rnspd\sth A\geq 0\} \]
for the subset of {\em doubly nonnegative} matrices. Trivially, any completely positive matrix $A$ is doubly nonnegative, and it is well-known \cite{BeSh} that the converse only holds for $n\leq 4$. Naturally, any rank-$k$ matrix $A\in\Rnspd$ can be decomposed as $A=CC^{\top}$, where $C$ is $n\times k$ but generally not nonnegative. This can easily be seen by using, for instance, the spectral theorem. If, additionally, there exists a $k\times k$ orthogonal matrix $Q$ that isometrically maps the columns of $C^{\top}$ in the nonnegative orthant $\RR^k_{\geq 0}$ of $\RR^k$, then the matrix $B^\top=QC^{\top}$ is nonnegative and
\be A=CC^{\top}=(CQ^{\top})(QC^{\top}) = BB^{\top}\ee
is a {\em completely positive factorization} of $A$. Note that even if such a $Q$ does not exist, $A$ may still be completely positive. To see this, let $C_m^{\top}$ be the $(k+m)\times n$ matrix that results when we add $m$ rows of zeros to $C^{\top}$. Then $A=C_mC_m^{\top}$ and there may exist an orthogonal matrix $Q$, now of size $(k+m)\times(k+m)$, such that $QC_m^{\top}$ is nonnegative. If $m$ is the smallest number of additional zero rows that is needed for such $Q$ to exist, then $k+m$ is the so-called {\em cp-rank} of $A$. As was shown in \cite{BaBe}, the maximum possible cp-rank $\phi(k)$ of any rank-$k$ completely positive matrix is bounded by $\half k(k+1)-1$, but the question which values can actually be attained is still an open problem. See also \cite{Sha}.

\begin{rem}\label{rem-1}{\rm Observe that $\phi(2)=2$. Indeed, it is easy to see that a subset $\UU\subset\RR^2$ can be rotated into $\RR^2_{\geq 0}$ if and only if the angle between each pair $u,v\in \UU$ is acute or right. For finite sets $\UU$ this shows that if $A\in\DD$ has rank $2$, then $A$ is completely positive with cp-rank $2$.}
\end{rem}

\subsection{Isometric embedding of vectors into the nonnegative octant of $\RR^3$}\label{sect-x}
The fact \cite{BaBe} that $\phi(3)=5$ shows that if $\UU\subset\RR^3$ cannot be simultaneously rotated into $\RR^3_{\geq 0}$,  after embedding $\UU$ isometrically in $\RR^4$ or $\RR^5$ it may be possible to rotate the embedded set into $\RR^4_{\geq 0}$ or $\RR^5_{\geq 0}$. An example of this counter-intuitive fact is given by the vectors  
\be\label{one} \small u_1=\left[\begin{array}{c} 2 \\ 0 \\ 0\end{array}\right], \hdrie
u_2=\left[\begin{array}{c} 0 \\ 2 \\ 0 \end{array}\right], \hdrie
u_3=\left[\begin{array}{c} 1 \\ 1 \\ \st\end{array}\right], \hdrie
u_4=\left[\begin{array}{c} 1 \\ 1 \\ -\st\end{array}\right]. \ee
Together with the origin, these are the vertices of a pyramid, that equals half a regular octahedron. See the left picture in Figure 1. Obviously, their mutual angles are nonobtuse.  
\begin{center} 
\begin{tikzpicture}[scale=1.1, every node/.style={scale=1.1}]

\draw[->] (0,0)--(3,0);
\draw[->] (0,0)--(0,1.5);
\draw[->] (0,0)--(-0.5,-2);

\node at (-0.3,0) {$o$};
\node at (3,-0.3) {$e_1$};
\node at (-0.3,1.5) {$e_2$}; 
\node at (-0.8,-2) {$e_3$};

\node at (2.2,0.2) {$u_1$};
\node at (-0.7,-1.2) {$u_2$};
\node at (1.2,1.1) {$u_3$};
\node at (1.3,-2) {$u_4$};

\draw[gray!40!white,fill=gray!40!white] (0,0)--(-0.33,-1.33)--(1,-1.9)--cycle;
\draw[fill=gray!20!white] (0,0)--(2,0)--(1,-1.9)--cycle;
\draw[fill=gray!40!white] (0,0)--(2,0)--(0.85,1)--cycle;
\draw[fill=gray!80!white] (0,0)--(0.85,1)--(-0.33,-1.33)--cycle;

\draw (-0.33,-1.33)--(1,-1.9)--(2,0);

\draw[fill=black] (0,0) circle [radius=0.05];
\draw[fill=black] (2,0) circle [radius=0.05];
\draw[fill=black] (-0.33,-1.33) circle [radius=0.05];
\draw[fill=black] (0.85,1) circle [radius=0.05];
\draw[fill=black] (1,-1.9) circle [radius=0.05];

\begin{scope}[shift={(7,-0.6)}]

\draw[->] (0,0)--(0,2);
\draw[->] (0,0)--(0,-1.5);
\draw[->] (0,0)--(-1.5,0);
\node at (0.3,0) {$o$};

\node at (0.3,1.5) {$e_2$};
\node at (-3,-0.3) {$e_3$};

\node[scale=0.9] at (-1.5,1.4) {$Pu_3$};
\node[scale=0.9] at (-1.5,-1.3) {$Pu_4$};
\node[scale=0.9] at (-2.3,0.3) {$Pu_2$};

\draw[fill=gray!80!white] (0,0)--(-2,0)--(-1,1.41)--cycle;
\draw[fill=gray!40!white] (0,0)--(-2,0)--(-1,-1.41)--cycle;

\draw[->] (0,0)--(-3,0);

\draw[fill=black] (0,0) circle [radius=0.05];
\draw[fill=black] (-1,1.41) circle [radius=0.05];
\draw[fill=black] (-1,-1.41) circle [radius=0.05];
\draw[fill=black] (-2,0) circle [radius=0.05];
\end{scope}

\begin{scope}[shift={(11.3,-0.6)}]

\draw[->] (0,0)--(0,2);
\draw[->] (0,0)--(0,-1);

\node at (0.3,1.5) {$e_2$};
\node at (-3,0.3) {$e_3$};
\node at (0.3,0) {$o$};

\begin{scope}[rotate=-35]
\node[scale=0.9] at (-1.7,1.1) {$RPu_3$};
\node[scale=0.9] at (-1.7,-1.7) {$RPu_4$};
\node[scale=0.9] at (-2.5,0.1) {$RPu_2$};

\draw[fill=gray!80!white] (0,0)--(-2,0)--(-1,1.41)--cycle;
\draw[fill=gray!40!white] (0,0)--(-2,0)--(-1,-1.41)--cycle;

\draw[fill=black] (0,0) circle [radius=0.05];
\draw[fill=black] (-1,1.41) circle [radius=0.05];
\draw[fill=black] (-1,-1.41) circle [radius=0.05];
\draw[fill=black] (-2,0) circle [radius=0.05];

\end{scope}
\draw[->] (0,0)--(-3,0);
\end{scope}

\end{tikzpicture} 
\end{center}
{\bf Figure 1. } The vectors $u_1,u_2,u_3,u_4$ (left); projections on the $(e_2,e_3)$ plane (middle) and the failed attempt to embed these into the nonnegative quadrant by a rotation $R$ (right).\\[2mm]
First observe that  
\be \small \half \left[\begin{array}{ccrr} \st & 0 & -1 & 1 \\ \st & 0 & 1 & -1\\ 0 & \st & 1 & 1 \\ 0 & \st & -1 & -1\end{array}\right]\left[\begin{array}{cccc} 2 & 0 & 1 & 1 \\ 0 & 2 & 1 & 1\\ 0 & 0 & \st & -\st \\ 0 & 0 & 0 & 0\end{array}\right] = \left[\begin{array}{cccc} \st & 0 & 0 & \st \\ \st & 0  & \st & 0\\ 0 & \st & \st & 0 \\ 0 & \st & 0 & \st\end{array}\right]. \ee
Thus, after embedding the vectors in (\ref{one}) isometrically in $\RR^4$, they can be isometrically mapped into $\RR^4_{\geq 0}$. There, of course, they still lie in a three-dimensional subspace.\\[2mm]
One particular way of proving that there exists no isometry $Q$ that maps $u_1,u_2,u_3,u_4$ into $\RR^3_{\geq 0}$, is to reduce the a priori {\em infinitely many} isometries $Q$ to be disqualified to only {\em finitely} many. The remaining ones can then one by one be inspected. In our above example, $u_1$ and $u_2$ are orthogonal, and thus so will be their images $Qu_1$ and $Qu_2$. These images can only be nonnegative if one of them is a positive multiple of one of the standard basis vectors $e_1,e_2,e_3$ of $\RR^3$, say $Qu_1=\|u_1\|e_1$, as is already the case. This reduces the problem to a two-dimensional rotation problem in the $(e_2,e_3)$-plane, which can be solved if and only if the orthogonal projections $Pu_2,Pu_3,Pu_4$ on that plane can be rotated into its nonnegative quadrant. But as explained in Remark \ref{rem-1}, this can be verified in a finite number of exact arithmetic operations. This is an example of a so-called {\em finiteness condition}. 

\begin{Pro}[Finiteness Condition I] \label{pro1} A point set in $\RR^3$ that contains two orthogonal vectors $u$ and $v$ can be rotated into the nonnegative octant if and only if it can be rotated into the nonnegative octant with at least one of the vectors $u$ and $v$ along a coordinate axis.
\end{Pro} 
In our specific example, $u_1$ is a positive multiple of $e_1$. The orthogonal projections of the remaining three vectors onto the $e_2,e_3$-plane, depicted in the middle of Figure 1, are
\be Pu_2=\left[\begin{array}{c} 2 \\ 0 \end{array}\right], \hdrie
Pu_3=\left[\begin{array}{c} 1 \\ \st\end{array}\right], \hdrie
Pu_4=\left[\begin{array}{c} 1 \\ -\st\end{array}\right]. \ee
The angle between $Pu_3$ and $Pu_4$ is obtuse. By symmetry, also the projections on the orthogonal complement of $e_2$, of which $u_2$ is a positive multiple, make an obtuse angle. Thus, $u_1,\dots,u_4$ cannot be isometrically embedded into $\RR^3_{\geq 0}$. See the right picture in Figure 1.

\subsection{Finiteness conditions and dimensional reduction}\label{sect-1.3}
For given $n\geq 2$, the problem of isometric embedding of a finite set into the nonnegative orthant of $\RR^n$ can be studied in its full generality as follows. Let $\UU=\{u_1,\dots,u_p\}\subset \RR^n$ be a finite point set. Obviously, $\UU$ can be isometrically embedded into $\RR^n_{\geq 0}$ if and only if there exists an orthonormal basis $\FF=\{f_1,\dots,f_n\}$ of $\RR^n$ such that all coordinates of each vector in $\UU$ with respect to the basis $\FF$ are nonnegative, or in other words, if and only if
\be\label{eq-1} f_i^\top u_j\geq 0 \hdrie \mbox{\rm for all $i\in\{1,\dots,n\}$ and all $j\in\{1,\dots,p\}$}. \ee
Writing $F$ for the matrix with columns $f_1,\dots,f_n$, and $U$ for the matrix with columns $u_1,\dots,u_p$, the condition in (\ref{eq-1}) is, of course, equivalent to $F^\top U\geq 0$ with $F^\top F=I$. 

\begin{rem}\label{rem-2}{\rm A necessary condition for the existence of $\FF$ is that $U^\top U\geq 0$, simply because inner products between vectors in $\RR^n_{\geq 0}$ are nonnegative. In $\RR^2$ this condition is also sufficient, as already stated in Remark \ref{rem-1}. In $\RR^n$ with $n\geq 3$ it is not.}
\end{rem}
Aiming for a recursive approach to solve the problem, observe that the set $\UU$ can be isometrically embedded into $\RR^n_{\geq 0}$ if and only if there exists a hyperplane $H$ with unit normal vector $g_1$ such that both of the following conditions, illustrated in Figure 2, are satisfied:\\[3mm]
(i) $\UU$ is a subset of the closed half space separated by $H$ in which $g_1$ lies as well,\\[2mm]
(ii) there exists an orthonormal basis $\GG=\{g_2,\dots,g_n\}$ of $H$ such that the set of orthogonal projections of the elements of $\UU$ onto $H$ have nonnegative coordinates with respect to $\GG$.\\[2mm]
Therefore, in theory, to solve an isometric embedding problem in $\RR^n$, it is sufficient to solve for each $g\in\mathbb{S}^{n-1}$ an isometric embedding problem in $\RR^{n-1}$. Now, a good finiteness condition is a practical condition that reduces the infinite number of vectors $g\in\mathbb{S}^{n-1}$ to be inspected to a finite number. In Section \ref{sect2}, we will formulate such a condition for $n=3$. 
\begin{center} 
\begin{tikzpicture}

\draw[gray!20!white, fill=gray!20!white] (4,0)--(8,0)--(9,1.5)--(5,1.5)--cycle;
\draw (0,0)--(8,0)--(10,3)--(2,3)--cycle;
\draw[->] (5,1.5)--(5,4.5);

\node at (5.4,4.3) {$g_1$};
\node at (3,1.5) {$H$};

\draw[thick,->] (5,1.5)--(7.5,2.3);
\draw[->] (5,1.5)--(7.5,1.5);
\draw[dotted](7.5,2.3)--(7.5,1.5);

\draw[thick,->] (5,1.5)--(4.7,4);
\draw[->] (5,1.5)--(4.7,0.2);
\draw[dotted] (4.7,4)--(4.7,0.2);

\draw[thick,->] (5,1.5)--(6.5,2.7);
\draw[->] (5,1.5)--(6.5,0.4);
\draw[dotted] (6.5,2.7)--(6.5,0.4);

\draw[thick,->] (5,1.5)--(8,1);
\draw[->] (5,1.5)--(8,0.4);
\draw[dotted] (8,1)--(8,0.4);

\draw[thick,->] (5,1.5)--(5.2,2.8);
\draw[->] (5,1.5)--(5.2,0.8);
\draw[dotted] (5.2,2.8)--(5.2,0.8);

\end{tikzpicture} 
\end{center}
{\bf Figure 2. } Dimensional reduction of the isometric embedding problem.\\[2mm]
Since the isometric embedding problem in $\RR^2$ can be explicitly solved by evaluating at most $4p$ inner products in $\RR^2$, this will lead to a finite exact algorithm to solve the problem in $\RR^3$.\\[2mm]
In accordance with Remark \ref{rem-2}, a necessary condition for the existence of $\GG$ in (ii) is that the mutual inner products between the projections of the elements of $\UU$ onto $H$ are nonnegative. Thus, the vector $g_1$ in (i) should be such that
\be\label{eq-2} U^\top U \geq  U^\top g_1g_1^\top U \geq 0. \ee 
Indeed, the right inequality in (\ref{eq-2}) shows that either $g_1$ or $-g_1$ satisfies (i), whereas the first implies that all inner products between vectors $(I-g_1g_1^\top)u_i$ and $(I-g_1g_1^\top)u_j$ are nonnegative.\\[3mm]
In $\RR^3$, the condition in (\ref{eq-2}) is also {\em sufficient} for the existence of $\GG$ in (ii).

\begin{Th} Let $\UU=\{u_1,\dots, u_p\}\subset \RR^3$. Then $\UU$ is isometrically embeddable in $\Rdp$ if and only if there exists a $g_1\in\mathbb{S}^2$ such that (\ref{eq-2}) holds. 
\end{Th}
{\bf Proof. } The question whether the projections of the vectors $u_1,\dots,u_p$ onto the two-dimensional orthogonal complement $H$ of $g_1$ can be isometrically embedded into a quadrant is equivalent with none of them making an obtuse angle, as stated in Remark \ref{rem-2}. \hfill $\Box$\\[2mm]
In Section \ref{sect2} we show that only {\em finitely many} $g_1$ need to be inspected in (\ref{eq-2}).
 
\section{Motivation}\label{sect1.2}
Our interest in isometric embedding of a set of vectors in $\RR^3$ into the nonnegative octant originates from the following problem. Consider the three-dimensional space $\PP^4_0(I)$ of polynomials of degree at most four that vanish at the boundary points $x=-1$ and $x=1$ of $I=[-1,1]$. The symmetric bilinear form
\be \langle p,q\rangle = \int_{-1}^1 p'(x)q'(x)dx\ee
defines an inner product on $\PP^4_0(I)$, the so-called $H^1_0(I)$-inner product. We wish to investigate if there exists an $\langle\cdot,\cdot\rangle$-orthonormal basis for $\PP_0^4(I)$ consisting of {\em nonnegative} polynomials. The existence of such a basis would imply a discrete maximum principle for certain finite element approximations of elliptic two-point boundary value problems \cite{VeSo}.

\subsection{Nonnegative $H^1_0(I)$-orthogonal polynomials}
Integration of the $L^2(I)$-orthogonal Legendre polynomials results in $\langle\cdot,\cdot\rangle$-orthonormal {\em Lobatto} polynomials $\phi_2,\phi_3,\phi_4\in\PP^4_0(I)$. These Lobatto polynomials are obviously not  nonnegative. See Figure 3 for a picture of their graphs and their explicit forms. Observe that they share a common factor $q(x)=(x+1)(x-1)$.
\begin{center}  
\begin{tikzpicture}

\begin{axis}[ticks=none, xmin=-1, xmax=1,ymin=-0.32,ymax=0.65, axis x line=middle, axis y line=middle,axis line style=-];
\addplot[domain=-1:1, samples=100]{0.234*(1-x)*(2.236*x-1)*(2.236*x+1)*(1+x)};
\addplot[domain=-1:1, samples=100]{0.612*(1+x)*(1-x)};
\addplot[domain=-1:1, samples=100]{0.791*(1-x)*x*(1+x)};
\end{axis}

\node at (1.5,5) {$\phi_2$};
\node at (1.5,0.5) {$\phi_4$};
\node at (1.5,3) {$\phi_3$};
\node[scale=0.9] at (-0.4,1.7) {$-1$};  
\node[scale=0.9] at (7,1.7) {$1$}; 

\node[scale=0.9] at (10.5,3.3) {\fbox{\begin{minipage}{7cm} $q(x)=(x+1)(x-1)$\\[3mm]
$\phi_2(x)=\frac14\sqrt{6}\,q(x)$\\[3mm]
$\phi_3(x) = \frac14 x\sqrt{10}\,q(x)$\\[3mm]
$\phi_4(x) = \frac{1}{16}(\sqrt{5}x+1)(\sqrt{5}x-1)\sqrt{14}\,q(x)$\end{minipage}}};

\end{tikzpicture}
\end{center}
{\bf Figure 3: } The $H^1_0(I)$-orthonogonal Lobatto polynomials $\phi_2,\phi_3,\phi_4$.\\[2mm]
Now, consider the curve $\CC\subset\RR^3$ defined as the image of 
\be \Phi: I \mapsto \RR^{3}: \hdrie x \rightarrow \left[\begin{array}{c} \phi_2(x)\\ \phi_3(x) \\ \phi_4(x)\end{array}\right]. \ee
Only if there exists an orthogonal transformation $Q$ such that $Q\CC\subset\Rdp$, the functions $\psi_2,\psi_3,\psi_4$ defined by
\be \Psi: I \rightarrow \RR^{3}: \hdrie x \mapsto \left[\begin{array}{c} \psi_2(x)\\ \psi_3(x) \\ \psi_4(x)\end{array}\right] = Q\Phi(x)\ee
constitute a nonnegative $\langle\cdot,\cdot\rangle$-orthonormal basis for $\PP_0^4(I)$. Indeed, $\Psi'(x)=Q\Phi'(x)$ and one can easily verify that the $\langle\cdot,\cdot\rangle$-orthonormality of $\psi_2,\psi_3,\psi_4$. The curve $\CC$ is displayed in the left picture of Figure 4. In the right picture of Figure 4, the canonical projection $\pi(\CC\setminus\{0\})$ of $\CC\setminus\{0\}$ onto the $2$-sphere $\mathbb{S}^2$ is depicted, where
\be \pi: \RR^3\setminus\{0\} \rightarrow \mathbb{S}^2: \hdrie x \mapsto \frac{x}{\|x\|}.\ee
To enhance the perspective in this picture, $\mathbb{S}^2$ is visualized by a number of randomly selected dots on its surface. Moreover, a number of points on the projected curve $\pi(\CC\setminus\{0\})$ are connected by a line with the origin of $\RR^3$. 
\begin{center}   
\includegraphics[height=6.0cm]{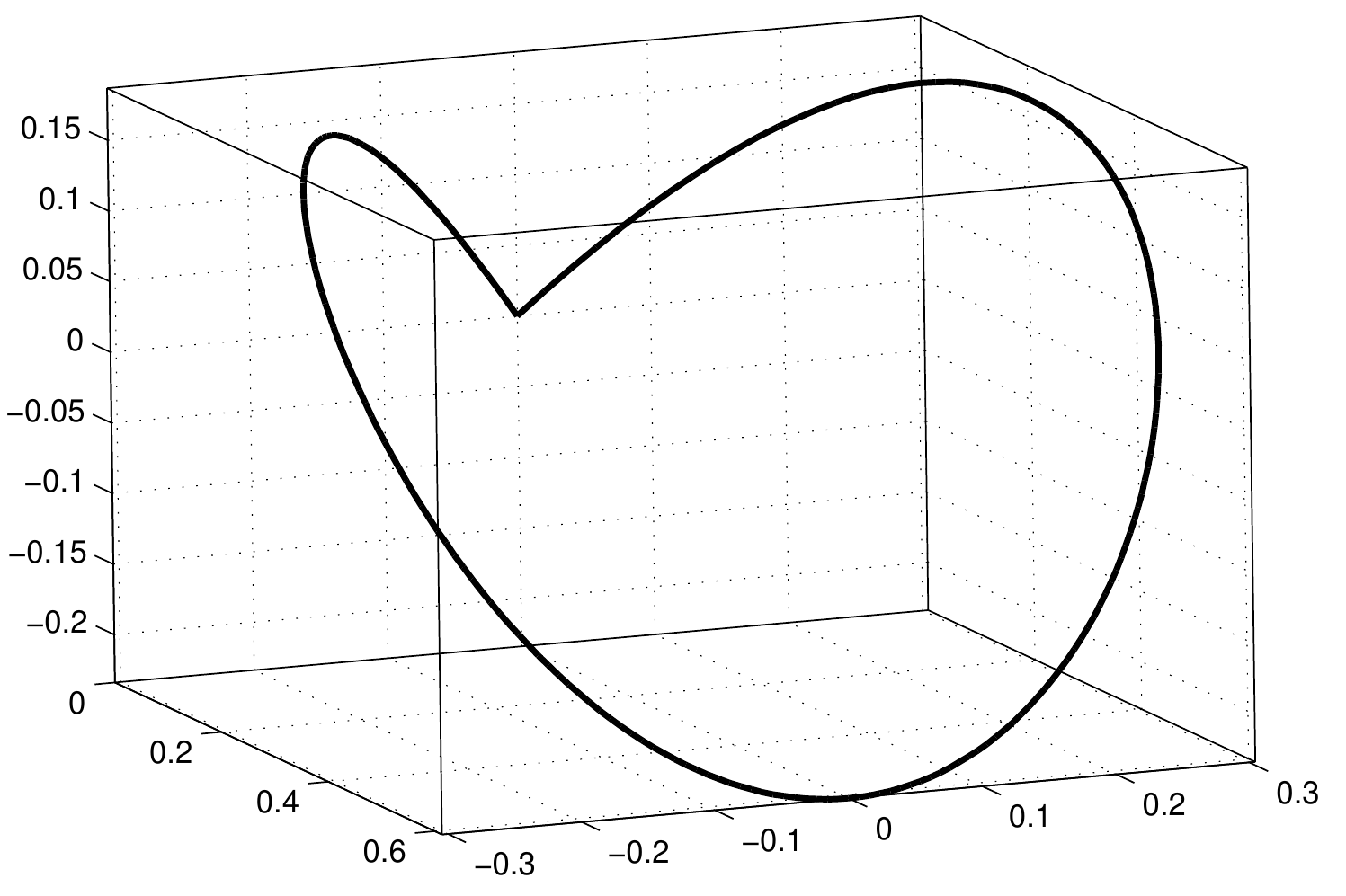}\hspace*{3mm}\includegraphics[height=6cm]{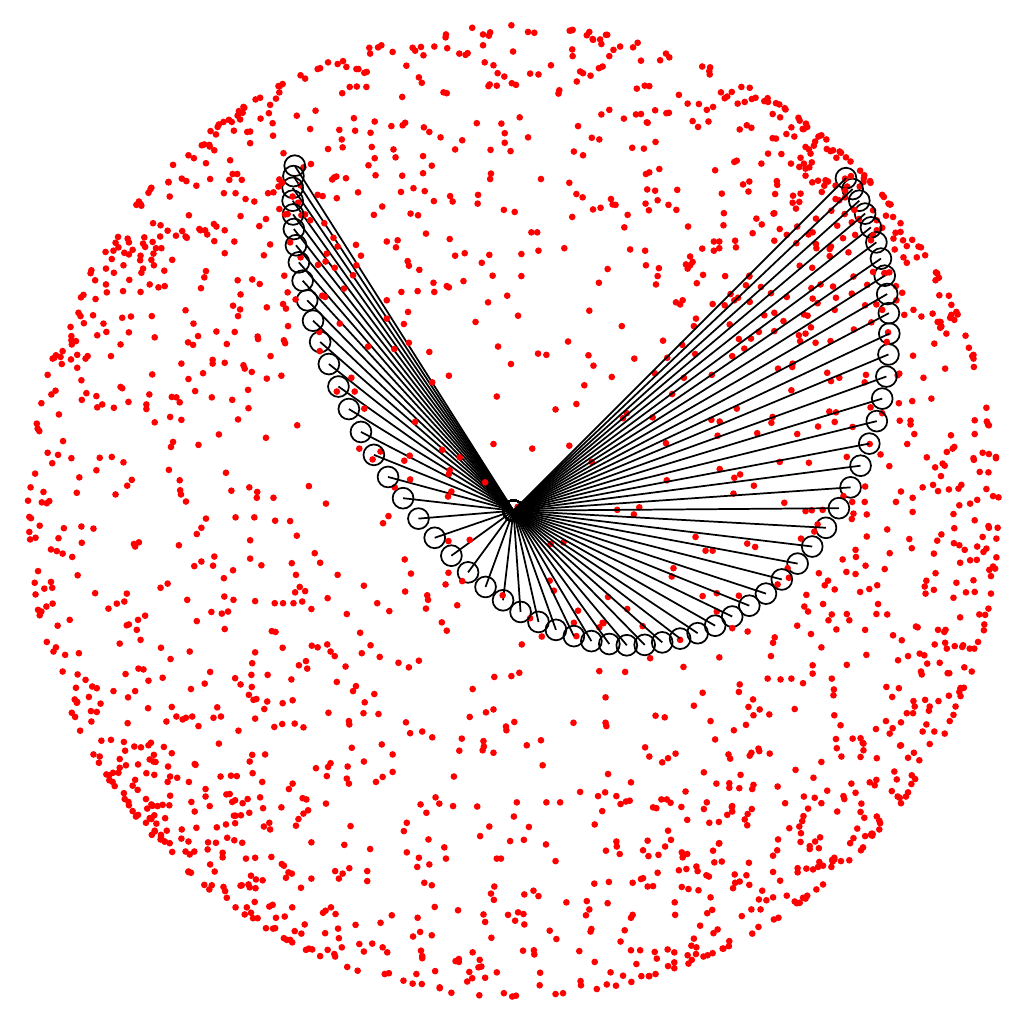}\\[3mm]
{\bf Figure 4. } Left: the graph $\CC\subset\RR^3$. Right: projected on the $2$-sphere $\mathbb{S}^2$.
\end{center} 
Clearly, for a given orthogonal transformation $Q$ we have that $QC\subset\Rdp$ if and only if  $Q\pi(\CC\setminus\{0\})\subset \Rdp\cap\mathbb{S}^2$. The mutual angles between points on $\CC$, are given by
\be \alpha: (-1,1)\times (-1,1)\rightarrow\RR: \hdrie (x,y)\mapsto  \arccos\left(\pi(\Phi(x))^\top\pi(\Phi(y))\right) \ee   
and visualized in Figure 5. None of them is obtuse, because it is easily verified that
\be \Phi(x)^\top\Phi(y)=\phi_2(x)\phi_2(y)+\phi_3(x)\phi_3(y)+\phi_4(x)\phi_4(y) \geq 0, \ee
and thus also
\be \pi(\Phi(x))^\top\pi(\Phi(y)) \geq 0.\ee
Hence, we cannot conclude from Remark \ref{rem-2} that the transformation $Q$ does {\em not} exist. As no mutual angle is right, Finiteness Condition I in Proposition \ref{pro1} can not be used either.
\begin{center} 
\includegraphics[height=5.5cm]{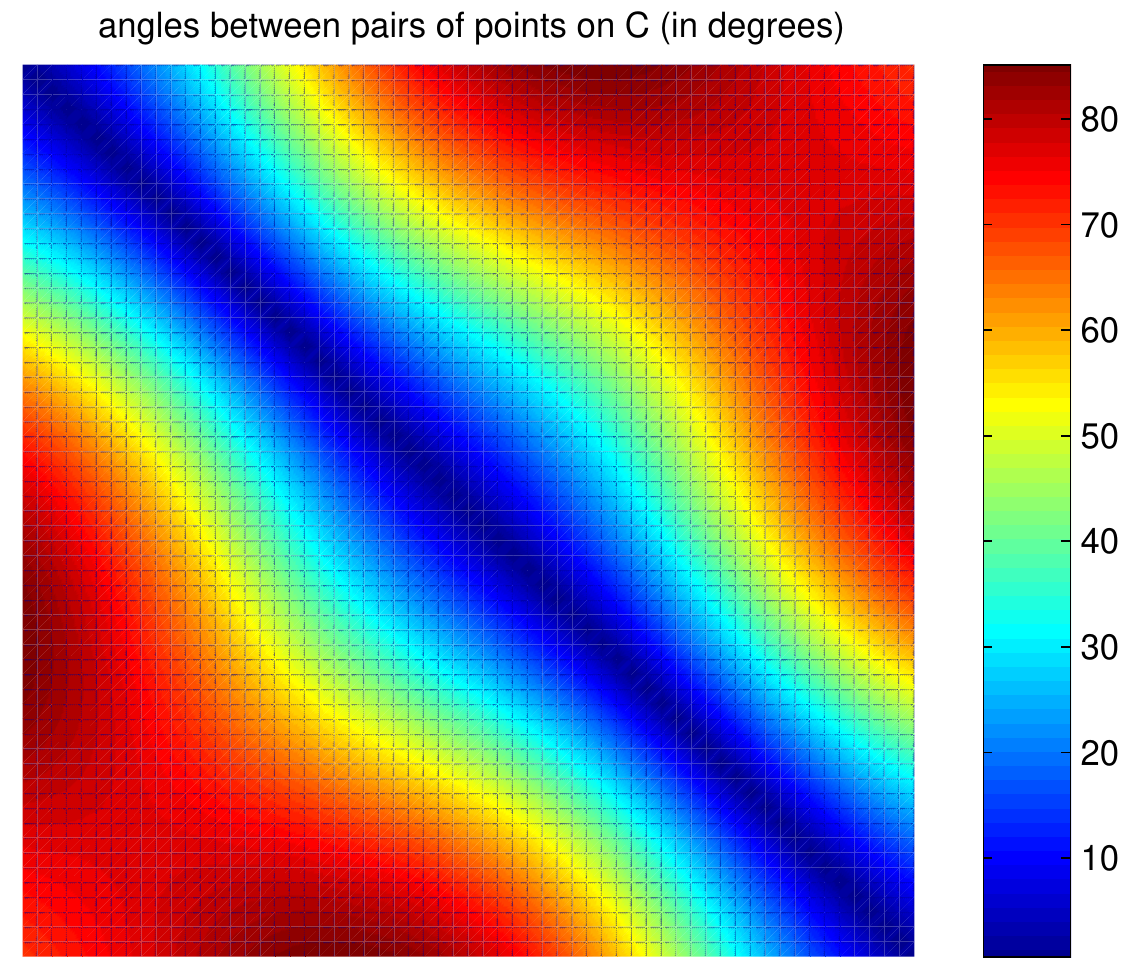}\hspace{3mm}\includegraphics[height=5.5cm]{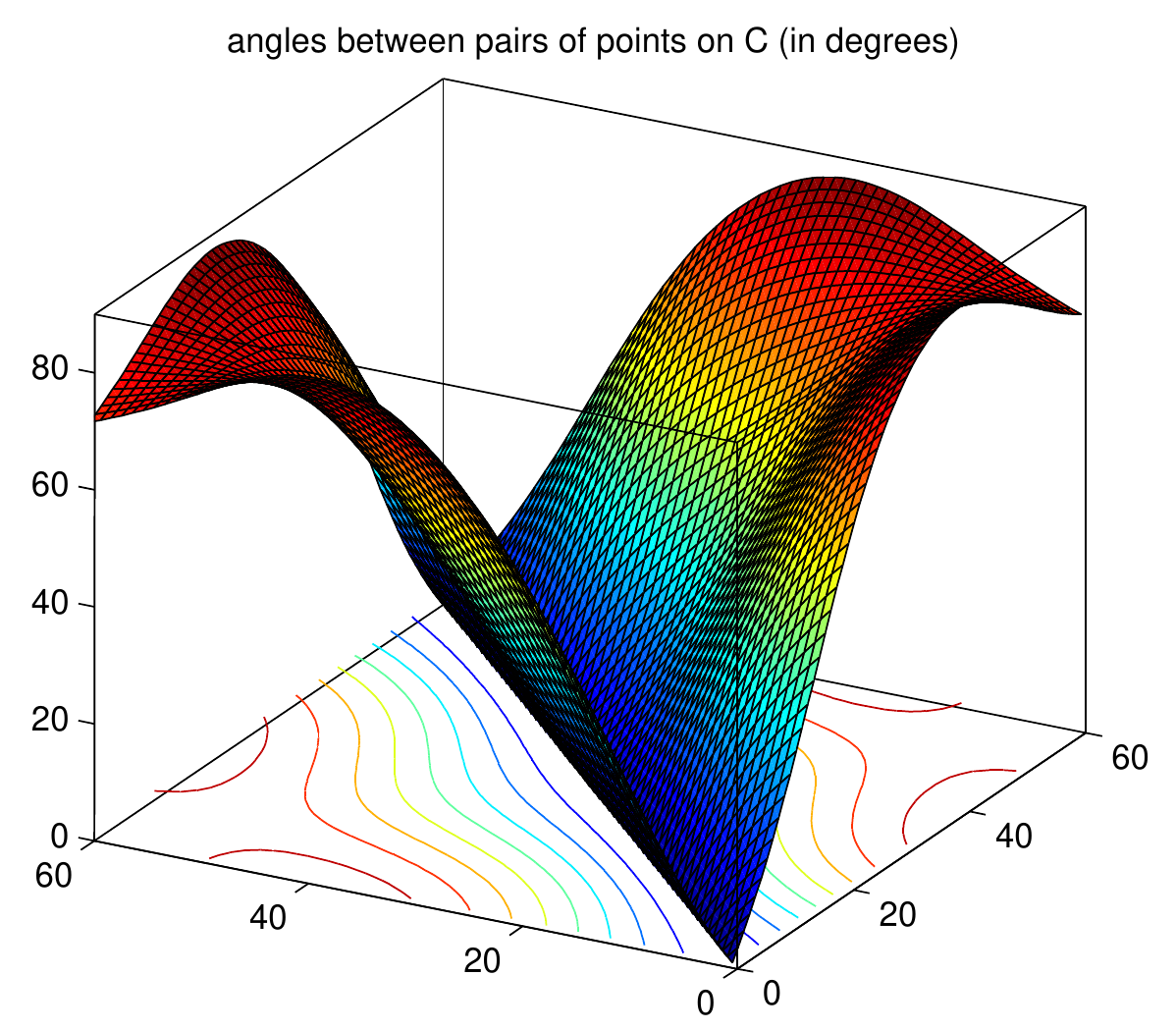}\\[2mm] 
{\bf Figure 5. } The angle $\alpha(x,y)$ between $\Phi(x)$ and $\Phi(y)$ does not reach $90^\circ$.
\end{center} 
Since also a direct analysis of the polynomials involved turned out to be quite tedious, we used an ad-hoc computer program to map the finite subsets     
\be \CC_\ell = \{ \Phi(x_j) \sth x_j = -1+\frac{j}{\ell}, \hdrie j=0,\dots,2\ell\}\ee
into $\Rdp$ for increasing values of $\ell$, using a discrete subset of the orthogonal transformations of $\RR^3$. For $\ell\leq 31$, thus with at most $63$ values of $\Phi$, we succeeded. The computed orthogonal transformation $\hat{Q}$ turned out to approximate the matrix $Q$ that is used in Figure 6 to transform $\Phi$ into $\Psi$ in two decimal places.
\begin{center} 
\begin{tikzpicture}

\begin{scope}
\begin{axis}[ticks=none, xmin=-1, xmax=1,ymin=-0.1,ymax=0.65,
      axis x line=middle, axis y line=middle,axis line style=-]
\addplot[domain=-1:1, samples=100]{0.40824829*0.612372*(1+x)*(1-x) +0.7071067*0.791*(1-x)*x*(1+x)   +0.5773502*0.234*(1-x)*(2.236*x-1)*(2.236*x+1)*(1+x)};
\addplot[domain=-1:1, samples=100]{0.40824829*0.612372*(1+x)*(1-x)  -0.7071067*0.791*(1-x)*x*(1+x)    +0.57735026*0.234*(1-x)*(2.236*x-1)*(2.236*x+1)*(1+x)};
\addplot[domain=-1:1, samples=100]{0.81649658*0.612372*(1+x)*(1-x)                -0.577350269*0.234*(1-x)*(2.236*x-1)*(2.236*x+1)*(1+x)};
\end{axis}
\end{scope}

\node at (2.5,5.5) {$\psi_2$};
\node at (6,4.5) {$\psi_4$};
\node at (1,4.5) {$\psi_3$};
\node[scale=0.9] at (-0.4,1) {$-1$};  
\node[scale=0.9] at (7,1) {$1$}; 

\node[scale=0.9] at (10.3,4) {$\left[\begin{array}{c}\psi_2 \\  \psi_3 \\ \psi_4\end{array}\right] = \left[\begin{array}{rrr}\frac{1}{3}\sqrt{6} & 0 & \frac{1}{3}\sqrt{3} \\ \frac{1}{6}\sqrt{6} & \frac{1}{2}\sqrt{2} & \frac{1}{3}\sqrt{3} \\ \frac{1}{6}\sqrt{6} & -\frac{1}{2}\sqrt{2} & \frac{1}{3}\sqrt{3}\end{array}\right]\left[\begin{array}{c}\phi_2 \\  \phi_3 \\ \phi_4\end{array}\right]$};

\node at (10.3,2) {\fbox{$\Psi = Q\Phi$}};
\end{tikzpicture}
\end{center}
{\bf Figure 6.} The isometrically transformed Lobatto polynomials $\psi_3,\psi_4$ are {\em not} nonnegative.\\[3mm]
Figure 6 suggests that we successfully transformed the polynomials $\phi_2,\phi_3,\phi_4$ into nonnegative $H^1_0$-orthonormal polynomials $\psi_2,\psi_3,\psi_4$. However, on closer inspection, $\psi_3$ and $\psi_4$ take negative values in the order of magnitude of $-5\times 10^{-4}$. Hence,  it is not clear if we should use an alternative discretization, or if the nonnegative orthonormal basis really does not exist.\\[3mm]
Even though $\CC$ consists of uncountably many points, we still may use the upcoming theory for finite sets in order to {\em disprove} the existence of the nonnegative orthonormal basis.

\begin{rem}\label{rem-2b}{\rm Note that $\phi_j=qr_j$ for all $j\in\{2,3,4\}$, where $q(x)=(x+1)(x-1)$ and $r_j$ is a polynomial of degree $j-2$. Since $q$ is nonnegative on $[-1,1]$, one may also study the problem of isometric embedding of the graph of
\[ R: I\rightarrow\RR^3: \hdrie x \mapsto \left[\begin{array}{r} r_2(x) \\ r_3(x) \\ r_4(x)\end{array}\right] \] 
into $\RR^3_{\geq 0}$. Observe that $R$ is simply a scaling of $\Phi$, hence the projection of its graph on $\mathbb{S}^2$ is the same as for the graph $\CC$ of $\Phi$, as depicted in Figure 4. Manipulations with polynomials in closed form are of course easier for $R$ than for $\Phi$, but they remain tedious.}
\end{rem}

\section{Finiteness conditions and containment problems}\label{sect2}
In this section, we formulate a generally applicable finiteness condition in the context of isometrically embedding a point set in $\RR^3$ into the nonnegative octant and prove its validity. For this, we generalize a so-called {\em containment problem} in plane geometry formulated in \cite{CrFaGu} to the corresponding result in spherical geometry. This problem was originally posed by Steinhaus in \cite{Ste}, solved by Post in \cite{Pos}, and generalized by Sullivan in \cite{Sul}

\begin{Th}[\cite{Sul}] Let $P$ be a polygon contained in a triangle $T$. Then $P$ also fits in $T$ with two of its vertices on the same edge of $T$. Moreover, the latter configuration can be realized using a continuous rigid transformation in which the polygon remains in $T$.
\end{Th}
It turns out that this result also holds for any {\em spherical} polygon $P$ contained in a right-angled equilateral {\em spherical} triangle $T$ on $\mathbb{S}^2$. Here, but also in \cite{Sul}, both $P$ and $T$ are supposed to be closed sets, and it is not necessary to assume that $P$ is (spherically) convex.

\subsection{Spherical polygon contained in a spherical triangle}\label{sect-2.1}
It seems nontrivial to modify the proof of Sullivan, which uses explicit calculations that involve orthogonal transformations that keep two vertices of a polygon on a pair of edges of the triangle. Fortunately, spherical triangles have certain properties that planar triangles do not have, and an easier alternative proof is available, as we shall see below. 
   
\begin{Def}{\rm For each $\alpha\in [0,\pi/2]$ we will write $\TT(\alpha)$ for the spherical triangle on $\sphere^2$ with vertices $t_1,t_2^\alpha,t_3$ given by
\[ t_1 = \left[\begin{array}{r} 1 \\ 0 \\ 0 \end{array}\right],\hdrie t_2^\alpha = \left[\begin{array}{c} \cos\alpha \\ \sin\alpha \\ 0 \end{array}\right],\und t_3 = \left[\begin{array}{r} 0 \\ 0 \\ 1 \end{array}\right]. \]
The edges of $\TT(\alpha)$ opposite $t_1,t_2^\alpha,t_3$, we denote by $\ell_1^\alpha,\ell_2,\ell_3^\alpha$, respectively. If $\alpha=\pi/2$ we omit $\alpha$ from the notation. In Figure 7 we depict $\TT$ and $\TT(\alpha)$ with $\alpha=\pi/6$.}
\end{Def}
\begin{center}   
\begin{tikzpicture}
\draw (3,1.5)--(1.5,0.75);
\draw[fill=gray!40!white] (0,0) to[out=90,in=210] (3,6)to[out=-110,in=90](2,-0.54) to[out=170,in=-20] (0,0);

\draw (0,0) to[out=-20,in=200] (6,0);
\draw (3,6) to[out=-30,in=90] (6,0);
\draw (2,-0.55)--(3,1.5);
\draw (3,1.5)--(6,0);
\draw (3,1.5)--(3,6);
\draw[fill] (0,0) circle [radius=0.05];
\draw[fill] (6,0) circle [radius=0.05];
\draw[fill] (3,6) circle [radius=0.05];
\draw[fill] (2,-0.54) circle [radius=0.05];
\draw[fill] (3,1.5) circle [radius=0.07];
\node[scale=1.2] at (-0.3,0) {$t_1$};
\node[scale=1.2] at (6.3,0) {$t_2$};
\node[scale=1.2] at (3.2,6.2) {$t_3$};
\node[scale=1.1] at (2.5,1) {$\alpha$};
\node at (1.3,2) {$\mathbb{T}(\alpha)$};
\node[scale=1.2] at (2,-1) {$t_2^\alpha$};
\node[scale=1.2] at (6,3) {$\ell_1$};
\node[scale=1.2] at (0,3) {$\ell_2$};
\node[scale=1.2] at (0.9,-0.7) {$\ell_3^\alpha$};
\node[scale=1.2] at (2.6,2.8) {$\ell_1^\alpha$};

\end{tikzpicture}
\end{center}
{\bf Figure 7. } The spherical triangle $\TT(\alpha)$ with $\alpha=\pi/6$ as subset of the spherical triangle $\TT$.\\[3mm]
Contrary to planar triangles, the spherical triangle $\TT$ has the property that the arc between a vertex $t_i$ and {\em any} point on $\ell_i$ has length $\pi/2$, and intersects $\ell_i$ orthogonally. Moreover, all three angles of $\TT$ are right. These properties help to prove Theorem \ref{th1} below.

\begin{Le}\label{lem-1} Let $\alpha\in [0,\pi/2]$. If a spherical polygon $P$ fits in $\TT(\alpha)$ with two of its vertices on the same edge of $\TT(\alpha)$, then $P$ also fits in $\TT$ with these two vertices on the same edge of $\TT$.
\end{Le}
{\bf Proof. } If the two vertices lie on $\ell_1^\alpha$, the only edge of $\TT(\alpha)$ that is generally not an edge of $\TT$, then a rotation about the $t_3$-axis over $\pi/2-\alpha$ maps $\ell_1^\alpha$ on $\ell_1$ while $P$ remains in $\TT$. \hfill $\Box$

\begin{Le}\label{lem-2} Let $\alpha\in [0,\pi/2]$. If a spherical polygon $P$ fits in $\TT(\alpha)$ with a vertex on a vertex of $\TT(\alpha)$, then $P$ also fits in $\TT$ with a vertex on a vertex of $\TT$.
\end{Le}
{\bf Proof. } The only vertex of $\TT(\alpha)$ that is generally not vertex of $\TT$ is $t_2^\alpha$. As in the previous lemma, a rotation about the $t_3$-axis over $\pi/2-\alpha$ maps $t_2^\alpha$ onto $t_2$ while $P$ remains in $\TT$. \hfill $\Box$
 
\begin{Th}\label{th1} Let $P$ be a spherical polygon contained in $\TT$. Then $P$ also fits in $\TT$ with two of its vertices on the same edge of $\TT$.
\end{Th} 
{\bf Proof. } Suppose that a spherical polygon $P$ fits into $\TT$. Then by compactness and continuity, the minimum
\be \beta = \min\left\{\alpha\in \left. \left[0,\half\pi\right] \hdrie \right|\hdrie \mbox{$P$ fits into $\TT(\alpha)$}\right\}\ee
exists. The actual configuration of $P$ within $\TT(\beta)$ does not need to be unique, but trivially, in each configuration, $P$ has at least one vertex on each edge of $\TT(\beta)$. Suppose that a vertex of $P$ lies on a vertex of $\TT(\beta)$. According to Lemma \ref{lem-2}, $P$ then fits into $\TT$ with a vertex on a vertex $p$ of $\TT$. Rotation about $p$ will move a second vertex of $P$ onto one of the two edges of $\TT$ that meet at $p$, and the theorem is proved. What remains is the case that exactly three vertices of $P$ lie on the boundary of $\TT(\beta)$, one on each edge. Denote the vertex of $P$ on $\ell_1^{\beta},\ell_2,\ell_3$ by $p_1^{\beta},p_2,p_3$, respectively. We will now construct a point $a\in\TT$ such that a rotation about the axis through the origin and $a$ moves $P$ into the interior of $\TT(\beta)$, contradicting the minimality of $\beta$. Indeed, let $a=\gamma_2\cap\gamma_3$, where $\gamma_j$ is the geodesic between $t_j$ and $p_j$ for $j\in\{2,3\}$. Note that these geodesics are orthogonal to $\ell_2$ and $\ell_3$, respectively. See Figure 8 for an illustration. 
\begin{center}   
\begin{tikzpicture}[scale=0.89]

\draw[fill=gray!20!white] (0,0) to[out=90,in=210] (3,6)to[out=-70,in=90](4,-0.54) to[out=185,in=-20](0,0);

\draw (0,0) to[out=90,in=210] (3,6)to[out=-70,in=90](4,-0.54);

\draw (0,0) to[out=-20,in=200] (6,0);
\draw (3,6) to[out=-30,in=90] (6,0);

\draw[fill] (0,0) circle [radius=0.05];
\draw[fill] (6,0) circle [radius=0.05];
\draw[fill] (3,6) circle [radius=0.05];
\draw[fill] (4,-0.54) circle [radius=0.05];

\node[scale=1.2] at (-0.3,0) {$t_1$};
\node[scale=1.2] at (6.3,0) {$t_2$};
\node[scale=1.2] at (3.2,6.2) {$t_3$};

\node[scale=1.2] at (4,-1) {$t_2^\beta$};

\node[scale=1.2] at (6,3) {$\ell_1$};
\node[scale=1.2] at (-0.1,2.3) {$p_2$};
\node[scale=1.2] at (1.7,-0.8) {$p_3$};
\node[scale=1.2] at (4.3,1.5) {$\ell_1^\beta$};

\node[scale=1.2] at (2,4) {$\gamma_3$};
\node[scale=1] at (1.3,3) {$\mathbb{T}(\beta)$};
\node[scale=1.2] at (3,1) {$\gamma_2$};
\node[scale=1.2] at (4.3,3.2) {$p_1^\beta$};

\draw (3,6) to[out=-113,in=91] (2,-0.54);
\draw (6,0) to[out=150,in=-10] (0.16,2);

\draw (2,1.65) to[out=45,in=210] (3.8,3);

\draw[thick,->] (3.8,3) to[out=110,in=-62] (3.3,4);
\draw[thick,->] (2,-0.54) to[out=-3,in=190] (3.2,-0.47);
\draw[thick,->] (0.16,2) to[out=-95,in=95] (0.16,1);

\draw[fill=white] (2,-0.54) circle [radius=0.07];
\draw[fill=white] (0.17,2) circle [radius=0.07];
\draw[fill=white] (3.8,3) circle [radius=0.07];

\draw[->] (2,1.65) circle [radius=0.7];

\draw[fill=white] (2,1.65) circle [radius=0.07];
\node[scale=1.2] at (1.7,1.4) {$a$};
\node[scale=1.2] at (3.6,2.5) {$\phi$};

\node[scale=1.2] at (3,2.8) {$\gamma$};

\end{tikzpicture}
\end{center} 
{\bf Figure 8. } Illustration of the main part of the proof of Theorem \ref{th1}.\\[3mm]
Next, consider also the geodesic $\gamma$ between $a$ and $p_1^{\beta}$ and write $\phi=\angle(a,p_1^\beta,t_2^\beta)$ for the angle that it makes with $\ell_1^{\beta}$. If $\phi\leq \pi/2$ then $P$ can be infinitesimally rotated about the axis through the origin and $a$ over a positive angle (in counter clockwise direction) such that all three vertices $p_1^{\beta},p_2,p_3$ move into $\TT(\beta)$ while $P$ remains in $\TT(\beta)$ and no vertices of $P$ are on the boundary of $\TT(\beta)$ anymore, obviously contradicting the minimality of $\beta$. Similarly, if $\phi\geq\pi/2$, then $P$ can be rotated over a negative angle, resulting in the same contradiction. Thus $P$ has two vertices on the same edge of $\TT(\beta)$, and Lemma \ref{lem-1} now finishes the proof. \hfill $\Box$\\[2mm]
As a corollary of Theorem \ref{th1}, we can now formulate a useful finiteness condition for the problem of isometric embedding of a finite point set in $\RR^3$ into its nonnegative octant.

\begin{Co}[Finiteness Condition II] Let $\UU=\{u_1,\dots,u_p\}\subset\RR^3$ be a set of $p$ pairwise linearly independent elements $(p\geq 2)$. Then $\UU$ can be isometrically embedded in $\Rdp$ if and only if this can be done with two of its elements having the same coordinate equal to zero.
\end{Co}
{\bf Proof. } Obviously, $\UU$ can be isometrically embedded in $\Rdp$ if and only if $\pi(\UU)$ can be isometrically embedded in $\TT$ on $\sphere^2$, and $\pi(\UU)$ can be isometrically embedded into $\TT$ if and only if the spherically convex hull $P$ of $\pi(\UU)$ can be isometrically embedded into $\TT$. Since $P$ is a spherical polygon, Theorem \ref{th1} proves the statement.\hfill $\Box$

\subsection{Efficient application of Finiteness Condition II}
The consequence of Finiteness Condition II in the context of Section \ref{sect-1.3} is, that  the set $\UU=\{u_1,\dots,u_p\}\subset \RR^3$ can be isometrically embedded in $\Rdp$ if and only if (\ref{eq-2}) holds for at least one of the $\half p(p-1)$ normal vectors $g_{ij}$ of the planes spanned by a pair $u_i,u_j$ with $i\not=j$.\\[3mm]
This number can be further reduced by the following observation. The convex hull $P$ of the $p$ vectors $\pi(u_1),\dots,\pi(u_p)$ on $\mathbb{S}^2$ is a convex spherical $m$-gon, with $m\leq p$. It is possible to list the $m$ vertices $\pi(u_{k_1}),\dots,\pi(u_{k_m})$ of $P$ in clockwise or counterclockwise order of traversal of the boundary of $P$ in a complexity of $\mathcal{O}(p\log\,m)$. See Section 1.1 of \cite{BeKrOvSc} for details. As a consequence, it is only needed to verify condition (\ref{eq-2}) for all $m$ normals $g_{k_j,k_{j+1}}$ to the planes spanned by the pairs $u_{k_j}$ and $u_{k_{j+1}}$ for $j\in\{1,\dots,m\}$, where $k_{m+1}=k_1$. Now, condition (\ref{eq-2}) can be verified by evaluating only $4p$ of the $\half p(p-1)$ mutual inner products between the projected vectors. An elegant way to do this is an inductive approach.
 
\begin{Pro} Suppose that from the set $\UU=\{u_1,\dots,u_{p-1}\}\subset\RR^2$ the vectors $u_1$ and $u_2$ make the largest mutual angle $\omega_{p-1}$. Assume that $\omega_{p-1}$ is not obtuse and that $u_1,u_2$ is positively oriented. Let $u_p\in\RR^2$ and compute
\be \alpha_1=u_1^\top u_p \und \alpha_2=u_2^\top u_p.\ee
Then the largest angle $\omega_p$ between all pairs from $\{u_1,\dots,u_p\}$ is non-obtuse if and only if $\alpha_1\geq 0$ and $\alpha_2\geq 0$. Moreover, at most one of the numbers $\beta_1,\beta_2$ defined by
\be \beta_1 = w_1^\top u_p \hdrie\mbox{\rm with}\hdrie w_1 = \left[\begin{array}{rr} 0 & -1 \\ 1 & 0\end{array}\right]u_2 \und
\beta_2 = w_2^\top u_p \hdrie\mbox{\rm with}\hdrie w_2 = \left[\begin{array}{rr} 0 & 1 \\ -1 & 0\end{array}\right]u_1\ee
can be negative, and\\[2mm]
$\bullet$ if $\beta_1<0$ then $\omega_p$ equals the angle between $u_1$ and $u_p$, and $u_1,u_p$ is positively oriented;\\[2mm]
$\bullet$ if $\beta_2<0$ then $\omega_p$ equals the angle between $u_p$ and $u_2$, and $u_p,u_2$ is positively oriented.\\[2mm]
If neither $\beta_1$ nor $\beta_2$ is negative, then $\omega_p$ equals the angle $\omega_{p-1}$ between $u_1$ and $u_2$.
\end{Pro}
The inductive approach is to keep track of a pair of angle maximizing vectors while considering increasingly more vectors from $\UU$. As soon as this angle would become obtuse, the process would be terminated. Otherwise, the final pair of angle maximizing vectors can be used to compute the rotation $Q$ into the nonnegative quadrant as follows.

\begin{Co} Suppose that from the set $\UU=\{u_1,\dots,u_p\}\subset\RR^2$ the vectors $u_1$ and $u_2$ make the largest mutual angle $\omega_p$, that $u_1,u_2$ are positively oriented, and that $\omega_p$ is not obtuse. Let $Q$ be the rotation matrix that maps $u_1$ onto $e_1$. Then $Q$ maps $\UU$ into $\RR^2_{\geq 0}$.
\end{Co}
Combining the above leads to the following complexity result for the implicit algorithm.

\begin{Th} Let $\UU=\{u_1,\dots,u_p\}\subset \RR^3$ and let $m\leq p$ be the number of vertices of the convex hull of $\pi(\UU)$ on $\mathbb{S}^2$. If there exists an orthonormal basis $\FF=\{f_1,f_2,f_3\}$ of $\RR^3$ such that the coordinates $f_i^\top u_j$ of all vectors in $\UU$ with respect to $\FF$ are nonnegative, the complexity of the algorithm that computes this basis is $\mathcal{O}(p^2\log m)$. The algorithm fails if $\FF$ does not exist.
\end{Th}

\subsection{Discussion of generalizations to higher dimensions}
The example in Section \ref{sect-x} shows that there exist finite subsets of $\RR^4_{\geq 0}$ that cannot be isometrically embedded in $\RR^4_{\geq 0}$ with three elements having the same coordinate equal to zero. As the four vectors lie in a hyperplane, this would imply that they all have the same coordinate equal to zero, and thus that they can be isometrically embedded in $\Rdp$. But in Section \ref{sect-x} we proved that this is impossible. The conclusion is that there is no straightforward generalization of Finiteness Condition II to dimension four. A strongly related observation is that even though the intersection of $\Rdp$ and a two-dimensional hyperplane always fits into $\RR^2_{\geq 0}$, the intersection of $\RR^4_{\geq 0}$ and a three-dimensional hyperplane does not necessarily fit into $\Rdp$. Alternatively, consider the intersection of the $3$-sphere $\mathbb{S}^3$ with $\RR^4_{\geq 0}$. This is a spherical tetrahedron $\mathbb{K}$ with only right dihedral angles. Intersecting it with a three-dimensional subspace yields a spherical triangle $T\subset\mathbb{K}$ that does not need to fit in one of the facets of $\mathbb{K}$. 
\begin{center}
\begin{tikzpicture}[scale=0.8]

\draw[fill=gray!20!white] (0,4)--(1.5,1)--(4,0)--(5.5,5)--cycle;

\draw (1.5,1)--(5.5,5);

\draw[gray!80!white,fill=gray!60!white,]  (0.75,2.5)--(2.75,0.5)--(4.75,2.5)--(2.75,4.5)--cycle;

\draw (4,0)--(0,4);
\draw (0,0)--(4,0)--(4,4)--(0,4)--cycle;
\draw (1.5,1)--(5.5,1)--(5.5,5)--(1.5,5)--cycle;
\draw (0,0)--(1.5,1);
\draw (4,0)--(5.5,1);
\draw (4,4)--(5.5,5);
\draw (0,4)--(1.5,5);

\draw[fill] (0.75,2.5) circle [radius=0.05];
\draw[fill] (2.75,0.5) circle [radius=0.05];
\draw[fill] (4.75,2.5) circle [radius=0.05];
\draw[fill] (2.75,4.5) circle [radius=0.05];

\begin{scope}[shift={(6,0)}]

\draw[fill=gray!20!white] (0,0)--(5.66,0)--(2.83,4.9)--cycle;

\draw[fill=gray!60!white] (1.41,0)--(4.24,0)--(4.24,2.83)--(1.41,2.83)--cycle;

\end{scope}
\end{tikzpicture}
\end{center}  
{\bf Figure 9. } The square fits into the regular tetrahedron, but not into a triangular facet.\\[3mm]
In fact, a similar statement is valid in Euclidean geometry. The intersection of a Euclidean tetrahedron with a plane does not necessarily fit into one of the facets of that tetrahedron. In Figure 9 we display a regular tetrahedron $K$ in a cube $C$ with edges of length $1$. Its intersection $S$ with a plane that halves four parallel edges of $C$ is a square with edges of length $\half\sqrt{2}$. This square does not fit into the equilateral triangle with edges of length $\sqrt{2}$. 
 
\begin{rem}{\rm Observe also that an infinitesimal perturbation of the square $S$ leads to a tetrahedron inside $K$ that does not fit into $K$ with three vertices on the same facet of $K$.}
\end{rem} 
 
\subsection{Solution of the motivational problem}
Application of Finiteness Condition II to the motivational problem described in Section \ref{sect2} leads to the following. Using the computer, we selected a subset of four candidate points on the curve $\CC$ of which we suspected that they cannot be isometrically embedded into $\Rdp$. These four points are the values of the function $R$ from Remark \ref{rem-2b} at the points $-1,-1/2,1/2,1$, which form the $3\times 4$ matrix $U$,
\be U = \left[\, R(-1) \sth R(-\half) \sth R(\half) \sth R(1)\,\right] =  \left[\begin{array}{ccc} \sqrt{6} & & \\ & \sqrt{10} & \\ & & \sqrt{14}\end{array}\right]\left[\begin{array}{rrcr} 1 & 1 & 1 & 1 \\ -1 & -\half & \half & 1 \\ 1 & \phantom{n}\frac{1}{16} & \frac{1}{16} & 1 \end{array}\right]. \ee
Note that the plane $y=0$ is a plane of symmetry of this set. See Figure 10 for an illustration. From this figure it is also clear that only for each of the three pairs $(R(\half),R(1))$ and $(R(-1),R(1))$ and $(R(-\half),R(\half))$ we need to verify if the projections of the four vectors on the plane spanned by this pair make an obtuse angle or not.
\begin{center}   
\begin{tikzpicture}[scale=1.3]

\draw[gray!40!white, fill=gray!40!white] (2.5,0.5)--(2.5,2.9)--(4,5.5)--(4,3.5)--cycle;

\draw[gray] (0.5,1)--(1.5,3);
\draw[gray] (1,2)--(4,2);
\draw[gray] (2.5,0.5)--(4,3.5);
\draw[gray] (1,2)--(1,5);

\draw[->] (1,2)--(4,5.5); 
\draw[->] (1,2)--(2.5,2.9);
\draw[->] (1,2)--(3.625,2.9);
\draw[->] (1,2)--(2.875,1.4);

\draw[dotted] (4,3.5)--(4,5.5);
\draw[dotted] (2.5,2.9)--(2.5,0.5);
\draw[dotted] (2.875,1.25)--(2.875,1.4);
\draw[dotted] (3.625,2.75)--(3.625,2.9);

\node[scale=0.8] at (0.8,2) {$0$};

\node[scale=0.8] at (4.4,5.5) {$\Phi(1)$};
\node[scale=0.8] at (3,3) {$\Phi(-1)$};
\node[scale=0.8] at (4.1,2.8) {$\Phi(\frac{1}{2})$};
\node[scale=0.8] at (3.4,1.3) {$\Phi(-\frac{1}{2})$};

\begin{scope}[shift={(6,1.5)}]

\draw[gray] (0,0)--(4,0);
\draw[gray] (2,0)--(2,3);

\draw[->] (2,0)--(0,3);
\draw[->] (2,0)--(4,3);
\draw[->] (2,0)--(1,0.16);
\draw[->] (2,0)--(3,0.16);

\node[scale=0.8] at (0.5,3) {$\Phi(-1)$};
\node[scale=0.8] at (1.2,0.4) {$\Phi(-\frac{1}{2})$};
\node[scale=0.8] at (2.9,0.4) {$\Phi(\frac{1}{2})$};
\node[scale=0.8] at (3.6,3) {$\Phi(1)$};

\end{scope}

\end{tikzpicture}
\end{center} 
{\bf Figure 10. } Left: The four vectors $R(-1), R(-\half), R(\half)$, and $R(1)$ in $\RR^3$, all lying in the (gray) plane $x=\sqrt{6}$. Right: view within the plane $x=\sqrt{6}$.\\[3mm]
An elementary calculation shows that each of three sets of projections has an obtuse angle among them. Hence, the four vectors cannot be isometrically embedded into $\Rdp$, and neither can the curve $\CC$. This confirms that the sought nonnegative orthonormal basis does not exist.

\begin{Th} There does not exist a nonnegative $H^1_0$-orthonormal basis for $\PP^4_0(I)$.
\end{Th}

\subsubsection*{Acknowledgments}
The authors acknowledge the support by Grant no.~P101/14-02067S of the Grant Agency of the Czech Republic and RVO 67985840. They also thank the anonymous referee for carefully studying the earlier versions of the manuscript and contributing to its improvement.

\end{document}